# GROMOV HYPERBOLIC JOHN IS QUASIHYPERBOLIC JOHN I

QINGSHAN ZHOU* AND SAMINATHAN PONNUSAMY

ABSTRACT. In this paper, we introduce a concept of quasihyperbolic John spaces and provide a necessary and sufficient condition for a space to be quasihyperbolic John. Using this criteria, we exhibit a simple proof to show that a John space with a Gromov hyperbolic quasihyperbolization is quasihyperbolic John, quantitatively. This answers affirmatively to an open question proposed by Heinonen (Rev. Math. Iber, 1989), which was studied by Gehring et al. (Math. Scand, 1989). As a tool, we study its connection between several geometric conditions.

## 1. Introduction and Main Results

This paper focuses on geometric properties of quasihyperbolic geodesics in John spaces. In what follows, $(X, d)$ is a locally compact, incomplete and rectifiably connected metric space and the identity map from $(X, d) \to (X, l)$ is continuous, where $l$ is the length metric of $X$ with respect to $d$. Following [1, 19, 24], for such a space, we call it minimally nice. Note that every proper domain of Euclidean spaces $\mathbb{R}^n$ is minimally nice as a space.

A minimally nice space $(X, d)$ is called *a-John* if there is a constant $a \geq 1$ such that every pair of two points $x, y$ in $X$ can be joined by an arc satisfying

$$\min\{\ell(\alpha[x, z]), \ell(\alpha[z, y])\} \leq a\, d(z),$$

where $d(z)$ is the distance between $z$ and the metric boundary $\partial X$ of $X$, $\alpha[x, z]$ and $\alpha[z, y]$ denote the two subarcs of $\alpha$ between $x$ and $z$, and $z$ and $y$, respectively. The arc $\alpha$ is called a *double a-cone* arc. The class of Euclidean John domains was first considered by John [22] in the study of elasticity theory.

In the literature, a minimally nice space $(X, d)$ is called *a-John with center $x_0$* if there is a constant $a \geq 1$ and a distinguished point $x_0 \in X$ such that for all $x \in X$, we may join $x$ to $x_0$ by a curve $\gamma$ in $X$ such that for all $z \in \gamma$,

$$\ell(\gamma[x, z]) \leq a\, d(z).$$

Note that the above definition ensures that $X$ is bounded. Indeed, if $X$ is bounded, then these two definitions are equivalent, see [15, 18, 27]. For a proof we refer to [29]. Our version of definition of John spaces works well even if $X$ is unbounded.

File: ZhouSamy5_2021_john_hyperbolicI.tex, printed: 2021-11-25, 12.06

2000 *Mathematics Subject Classification.* Primary: 30C65, 30L10; Secondary: 30F45.

*Key words and phrases.* John spaces, quasihyperbolic geodesic, Gromov hyperbolic spaces, uniform spaces.

* Corresponding author.

The research was partly supported by NSF of China (No. 11901090).

The quasihyperbolic metric was introduced by Gehring and Palka [14] in Euclidean setting and we refer to [13] for more geometric properties. In [26], Martio and Sarvas introduced the notion of uniform domains, i.e., every pair of points can be connected by a double cone arc satisfying the quasiconvexity condition. The importance of the quasihyperbolic metric and uniform domains in quasiconformal mapping theory is well understood, see for example [6, 8, 9, 11, 15, 17, 19, 21].

In [3], Bonk et al. investigated negative curvature of uniform metric spaces and demonstrated many phenomena in function theory from the point of view of Gromov hyperbolicity of the quasihyperbolic metric. It was shown in [3, Theorem 2.10] that every quasihyperbolic geodesic in a uniform space is a uniform arc. Note the existence of a quasihyperbolic geodesic between any pair of points in a minimally space follows from [3, Proposition 2.8]. Since the uniformity of spaces implies the John property, it raises the natural problem of determining sufficient and/or necessary conditions (in a John space) for quasihyperbolic geodesics to be double cone arcs.

This problem has been studied by several authors, for background and more information see [29]. For example, in the case of plane subdomains, Gehring et al. [12] proved that simple connectedness is a sufficient condition. In the case of dimensions $n \geq 3$, Gromov hyperbolicity of the quasihyperbolic metric is also sufficient to assert that quasihyperbolic geodesics are double cone arcs. Using conformal modulus of path families and Ahlfors $n$-regularity of $n$-Lebesgue measure, this result was established in [3, Proposition 7.12] via the fact that John domains of $\mathbb{R}^n$ satisfy linearly locally connected condition with respect to the length metric. For a new proof of this result by using geometric characterization of Gromov hyperbolicity established in [23, 24], see [17].

In this paper, we focus on this problem and introduce the following concept.

**Definition 1.1.** Let $a \geq 1$. A minimally nice space $(X, d)$ is called *quasihyperbolic a-John*, if every quasihyperbolic geodesic $\alpha$ in $X$ is a double $a$-cone arc. Moreover, $X$ is called *quasihyperbolic John* if there is a constant $a \geq 1$ such that $X$ is quasihyperbolic $a$-John.

Obviously, quasihyperbolic John spaces are John. There are John spaces which are not quasihyperbolic John, see [12, Examples 5]. It follows from [3, Theorem 2.10] that every uniform space is quasihyperbolic John. Indeed, we shall see that the class of quasihyperbolic John spaces is very wide.

In the following, we first present a complete solution to the above problem by providing a necessary and sufficient condition for spaces to be quasihyperbolic John.

**Theorem 1.2.** *A minimally nice space $X$ is quasihyperbolic $a$-John if and only if there is a constant $A > 0$ such that for each quasihyperbolic geodesic $[x,y]_k$ satisfying $d(u) \leq 2\min\{d(x), d(y)\}$ for all $u \in [x,y]_k$, we have $k(x,y) \leq A$, where $k$ is the quasihyperbolic metric of $X$. The constants $a$ and $A$ depend only on each other.*

Next we consider whether Theorem 1.2 can tell us that Gromov hyperbolic John spaces are quasihyperbolic John. Note that, in a minimally nice space $(X, d)$, $X$ is called Gromov $\delta$-hyperbolic if it has a *Gromov hyperbolic quasihyperbolization*. That is to say, $(X, k)$ is $\delta$-hyperbolic for some constant $\delta \geq 0$, where $k$ is the quasihyperbolic metric of $X$. The class of Gromov hyperbolic John spaces includes uniform domains and inner uniform domains of $\mathbb{R}^n$, simply connected John domains in the plane, and John domains of $\mathbb{R}^n$ which are Gromov $\delta$-hyperbolic in the quasihyperbolic metric, etc.

In particular, it was asked by Heinonen in [18, Question 2] whether a John domain of $\mathbb{R}^n$ which quasiconformally equivalent to the unit ball or a uniform domain is quasihyperbolic John. Recently, the authors in [29] obtained a dimension-free answer to this question by showing that every quasihyperbolic geodesic in a John space admitting a roughly starlike Gromov hyperbolic quasihyperbolization is a double cone arc. The approach in [29] is elementary and makes use of the uniformization theory of Gromov hyperbolic spaces established in [3]. With the aid of Theorem 1.2, we obtain the following:

**Theorem 1.3.** *Let $a \geq 1$ and $\delta \geq 0$. Every quasihyperbolic geodesic in a minimally nice $a$-John Gromov $\delta$-hyperbolic space $X$ is a double $M$-cone arc with $M = M(a, \delta)$. In particular, a John space with a Gromov hyperbolic quasihyperbolization is quasihyperbolic John.*

It follows from Theorem 1.3 that quasihyperbolic geodesics in John hyperbolic spaces are double cone arcs, quantitatively. Moreover, Theorem 1.3 is an improvement of [3, Proposition 7.12], [12, Theorem 4.1] and [17, Remark 3.10]. We remark that the related results of [3, 12, 17] were considered in bounded domains of Euclidean spaces $\mathbb{R}^n$ and the parameters depended on $n$. Indeed, their proofs may not hold without assuming the spaces to be Ahlfors regular. This is because we do not know whether the spaces/domains satisfy the ball separation condition in this setting.

On the other hand, the rough starlikeness is not required in Theorem 1.3 and so it gives an improvement of [29, Theorem 1.1]. Indeed, our proof is not only more direct but is also simpler than that of [29, Theorem 1.1]. An important observation is that every double cone arc is a union of two quasihyperbolic quasigeodesics. Using this observation and the stability property of quasigeodesics in Gromov hyperbolic spaces, one can easily verify the criteria for quasihyperbolic John spaces given in Theorem 1.2. See Lemma 3.2.

In view of the above discussions, Theorem 1.3 asserts that simply connected plane John domains and Gromov hyperbolic John domains of $\mathbb{R}^n$ are quasihyperbolic John. There are many applications of these domains in quasiconformal mapping theory and potential analysis. For example, Ghamsari et al. [15] proved that quasisymmetric (resp. bilipschitz) maps between boundaries of John disks can be extended to a quasiconformal (resp. bilipschitz) map of the extended plane. In [17], Guo studied the uniform continuity of quasiconformal mappings onto Gromov hyperbolic generalized $\varphi$-John domains. The crucial ingredient in their proofs is based on the fact that quasihyperbolic geodesics in a Gromov hyperbolic John (or $\varphi$-John) domain are inner ($\varphi$-inner) uniform curves. Recently, Chen and Ponnusamy [11] established certain relationship between $K$-quasiconformal harmonic mappings and John disks.

Moreover, we study the connection between quasihyperbolic John spaces and other geometric conditions, such as linearly locally connected condition, the Gehring-Hayman condition and the ball separation condition. We begin with some definitions. Following [1] and [5], we say that $(X, d)$ is a *GHS space*, if it is a minimally nice length space and it satisfies both the Gehring-Hayman condition and the ball separation condition, see Definitions 3.4 and 3.5. A minimally nice space $(X, d)$ is $c$-$LLC_2$, if $c \geq 1$ and points in $X \setminus \overline{B}(x, r)$ can be joined by a curve in $X \setminus \overline{B}(x, r/c)$. In the following, we demonstrate that a GHS space satisfying the John or $LLC_2$ condition is quasihyperbolic John.

**Theorem 1.4.** *Let $(X, d)$ be a minimally nice length space. If $X$ is a GHS space satisfying the John or $LLC_2$ condition, then $X$ is quasihyperbolic John.*

As applications of the above results, we investigate geometric characterizations of uniform metric spaces. Here is an immediate corollary to Theorem 1.4.

**Corollary 1.5.** *A minimally nice length space is uniform if and only if it is John and GHS.*

As a corollary to Theorem 1.3, we investigate the connection between uniform spaces and quasihyperbolic John spaces under the assumption of Ahlfors regularity.

**Corollary 1.6.** *Let $(X, d, \mu)$ be a minimally nice, Ahlfors regular metric measure space. Then $X$ is uniform if and only if it is quasiconvex, John and has a roughly starlike Gromov hyperbolic quasihyperbolization.*

Using a similar argument as in the proof of [17, Proposition 3.7], one may easily get Corollary 1.6 provided the underlying space is locally externally connected (cf. [5]). Since not all of uniform spaces are locally externally connected, this condition is unnecessary. Under the absence of this extra condition, Corollary 1.6 follows from Theorem 1.3 and the Gehring-Hayman condition (cf. [24]).

In [1], Balogh and Buckley proved the equivalence of three different geometric properties of metric measure spaces, that is, Gromov hyperbolicity of the quasihyperbolic metric, the slice condition (for precise definition see [1]) and a combination of the Gehring-Hayman and ball separation conditions. So one may characterize uniform spaces in terms of the slice condition, for related discussions see for example [5, 6, 9, 10] and the references therein. Along this way, we now present a second application of Theorem 1.3. We also remark that the slice condition is unnecessary for uniformity with the absence of the locally externally connectedness, see [10].

**Corollary 1.7.** *A quasiconvex, John minimally nice metric space $(X, d)$ satisfying the slice condition is uniform.*

This paper is organized as follows. Section 2 contains notations and the basic definitions. After establishing certain useful lemmas, we prove Theorems 1.2, 1.3 and 1.4 in Section 3. Section 4 is devoted to the proofs of Corollaries 1.6 and 1.7.

## 2. Preliminaries

Let $(X, d)$ denote a metric space with its metric completion by $\bar{X}$ and its metric boundary by $\partial X = \bar{X} \setminus X$. The space $X$ is incomplete if $\partial X$ is non-empty. For all $x \in X$, $d(x) = \text{dist}(x, \partial X)$ if $\partial X \neq \emptyset$. The open (resp. closed) metric ball with center $x \in X$ and radius $r > 0$ is denoted by

$$B(x, r) = \{z \in X : d(z, x) < r\} \text{ (resp. } \overline{B}(x, r) = \{z \in X : d(z, x) \leq r\}),$$

and the metric sphere by $S(x, r) = \{z \in X : d(z, x) = r\}$.

A *curve* in $X$ is a continuous map $\gamma : I \to X$ from an interval $I \subset \mathbb{R}$ to $X$. If $\gamma$ is an embedding of $I$, it is also called an *arc*. We also denote the image set $\gamma(I)$ of $\gamma$ by $\gamma$. The *length* $\ell(\gamma)$ of $\gamma$ with respect to the metric $d$ is defined in an obvious way. Here the parameter interval $I$ is allowed to be open or half-open. Suppose $\gamma$ is a curve in $X$ with endpoints $x$ and $y$. We say that $z$ is a *midpoint* of $\gamma$ if $\ell(\gamma[x, z]) = \ell(\gamma[y, z])$.

A *geodesic* arc $\gamma$ joining $x \in X$ to $y \in X$ is a continuous map $\gamma$ from an interval $I = [0, l] \subset \mathbb{R}$ into $X$ such that $\gamma(0) = x$, $\gamma(l) = y$ and

$$d(\gamma(t), \gamma(t')) = |t - t'| \text{ for all } t, t' \in I.$$

A metric space $X$ is said to be *geodesic* if every pair of points can be joined by a geodesic arc. We use $[x, y]$ to denote a geodesic between two points $x$ and $y$ in $X$. Let $\lambda \geq 1$ and $\mu \geq 0$. A $(\lambda, \mu)$-*quasigeodesic* curve in $X$ is a $(\lambda, \mu)$-quasi-isometric embedding $\gamma : I \to X$, $I \subset \mathbb{R}$. More explicitly,

$$\lambda^{-1}|t - t'| - \mu \leq d(\gamma(t), \gamma(t')) \leq \lambda|t - t'| + \mu \text{ for all } t, t' \in I.$$

A metric space $(X, d)$ is called *rectifiably connected* if every pair of points in $X$ can be joined with a curve $\gamma$ in $X$ with $\ell(\gamma) < \infty$. In a minimally nice space $(X, d)$, the *quasihyperbolic metric* $k$ in $X$ is defined by

$$k(x, y) = \inf\left\{\int_\gamma \frac{ds}{d(z)}\right\},$$

where the infimum is taken over all rectifiable curves $\gamma$ in $X$ with endpoints $x$ and $y$, and $ds$ denotes the arc-length element with respect to the metric $d$. We remark that the resulting space $(X, k)$ is complete, proper and geodesic (cf. [3]). Moreover, we need the following well-known inequalities (cf. [3, 5]): for all $x, y \in X$,

$$(2.1) \qquad \left|\log \frac{d(x)}{d(y)}\right| \leq k(x, y)$$

and

$$(2.2) \qquad \log\left(1 + \frac{\ell(\gamma)}{\min\{d(x), d(y)\}}\right) \leq \ell_k(\gamma),$$

where $\gamma$ is a curve in $X$ with endpoints $x, y \in \gamma$ and $\ell_k(\gamma)$ is the quasihyperbolic length of $\gamma$.

Let $\delta \geq 0$. A geodesic space $(X, d)$ is called (Gromov) $\delta$-*hyperbolic* if for each triples of geodesics arcs $[x, y], [y, z], [z, x]$ in $(X, d)$, every point in $[x, y]$ is within distance $\delta$ from $[y, z] \cup [z, x]$. Let $X$ be a proper, geodesic $\delta$-hyperbolic space, and let $w \in X$ and $K \geq 0$. We say that $X$ is $K$-*roughly starlike* with respect to $w$ if for each $x \in X$ there is a geodesic ray $\gamma = [w, \xi]$ connecting $w$ to $\xi$ such that $\text{dist}(x, \gamma) \leq K$, where $\xi$ is a point on the Gromov boundary of the $\delta$-hyperbolic space $X$.

We now recall the stability property of quasigeodesics in a $\delta$-hyperbolic space.

**Theorem 2.1.** ([4, Chapter III.H Theorem 1.7]) *For all $\delta \geq 0, \lambda \geq 1, \mu \geq 0$, there is a number $R = R(\delta, \lambda, \mu)$ with the following property: If $X$ is a $\delta$-hyperbolic space, and if $\gamma$ is a $(\lambda, \mu)$-quasigeodesic in $X$ and $[p, q]$ is a geodesic joining the endpoints of $\gamma$, then the Hausdorff distance between $[p, q]$ and $\gamma$ is less than $R$.*

Let $a \geq 1$. A metric space $X$ is called *a-quasiconvex* if each pair of points $x, y \in X$ can be joined by an $a$-quasiconvex curve $\alpha$, that is, $\ell(\alpha) \leq a\, d(x, y)$. A minimally nice space $X$ is called *a-uniform* if each pair of points in $X$ can be joined by a double $a$-cone and $a$-quasiconvex arc.

## 3. Proofs of Main results

### 3.1. A criteria for quasihyperbolic John spaces.
In this subsection, we study geometric properties of quasihyperbolic geodesics in quasihyperbolic John spaces. The aim of this part is to show Theorem 1.2 which provides a criteria for quasihyperbolic John spaces. Before the proof, we need an auxiliary lemma which asserts that every double cone arc is a union of two quasihyperbolic quasigeodesics.

**Lemma 3.1.** *Let $(X, d)$ be a minimally nice $a$-John space, and let $\gamma$ be a double $a$-cone arc with endpoints $x_1$ and $x_2$ in $X$, and let $x_0$ be a midpoint of $\gamma$. Then for any $y \in \gamma[x_1, x_0]$ and $z \in \gamma[y, x_0]$, we have*

$$(3.1) \qquad k(y, z) \leq \ell_k(\gamma[y, z]) \leq 3a \log\left(1 + \frac{a d(z)}{d(y)}\right) \leq 3a k(y, z) + 3a \log(3a).$$

*Moreover, the subarcs $\gamma[x_1, x_0]$ and $\gamma[x_0, x_2]$ are both quasihyperbolic $(\lambda, \mu)$-quasigeodesics with $\lambda = 3a$ and $\mu = 3a \log(3a)$.*

*Proof.* We only need to verify (3.1), because from which it follows that the subarc $\gamma[x_1, x_0]$ is evidently a quasihyperbolic quasigeodesic; and by using a symmetric argument, one easily see that the other arc $\gamma[x_0, x_2]$ shares the same property.

Towards this end, since $\gamma$ is a double $a$-cone arc and since $x_0$ is a midpoint of $\gamma$, for any $u \in \gamma[y, z]$, we have

$$(3.2) \qquad \ell(\gamma[y, u]) \leq a d(u).$$

Also, we claim that

$$(3.3) \qquad d(y) \leq 2a d(u).$$

Indeed, if $u \in B(y, d(y)/2)$, then we get

$$d(u) \geq d(y) - d(y, u) \geq d(y)/2.$$

Otherwise, by (3.2) we have

$$a d(u) \geq \ell(\gamma[y, u]) \geq d(y, u) \geq d(y)/2,$$

as required.

Now, by (3.2) and (3.3), we obtain

$$\ell(\gamma[y, u]) + d(y) \leq 3a d(u).$$

Therefore,

$$\begin{aligned}
k(y,z) &\leq \ell_k(\gamma[y,z]) = \int_{\gamma[y,z]} \frac{|du|}{d(u)} \\
&\leq \int_0^{\ell(\gamma[y,z])} \frac{3a\,dt}{t+d(y)} = 3a\log\left(1+\frac{\ell(\gamma[y,z])}{d(y)}\right) \\
&\leq 3a\log\left(1+\frac{ad(z)}{d(y)}\right) \\
&\leq 3a\log\frac{d(z)}{d(y)} + 3a\log(3a) \\
&\leq 3ak(y,z) + 3a\log(3a)
\end{aligned}$$

and the lemma follows. $\square$

**Proof of Theorem 1.2.** We first prove the necessity. Suppose that $X$ is quasi-hyperbolic $a$-John with constant $a > 0$. Fix $x, y \in X$ and let $[x,y]_k = \alpha$ be a quasi-hyperbolic geodesic connecting $x$ and $y$, and satisfying $d(u) \leq 2\min\{d(x), d(y)\}$ for all $u \in [x,y]_k$. Let $x_0$ be the midpoint of $\alpha$. By Lemma 3.1, we have

$$k(x,x_0) \leq 3a\log\left(1+\frac{ad(x_0)}{d(x)}\right) \leq 3a\log(1+2a)$$

and similarly, $k(x_0,y) \leq 3a\log(1+2a)$. These two estimates imply that $k(x,y) \leq 6a\log(1+2a) =: A$, as desired.

We are thus left to show the sufficiency. Let $[p,q]_k$ be a quasihyperbolic geodesic in $X$ with endpoints $p$ and $q$. We show that $[p,q]_k$ is a double $M$-cone arc with $M = M(A)$. Let $w_0 \in [p,q]_k$ be a point such that

$$d(w_0) = \max\{d(z): \ z \in [p,q]_k\} := T.$$

There is a unique nonnegative integer $n$ such that

$$2^n d(p) \leq T < 2^{n+1} d(p).$$

For each $i = 0, \ldots, n$, let $p_i$ be the first point on $[p,q]_k$ with

(3.4) $$d(p_i) = 2^i d(p),$$

when traveling from $p$ to $q$.

Similarly, we define $q_j$ to be the first point on $[p,q]_k$ with

$$d(q_j) = 2^j d(q)$$

for $j = 0, \ldots, m$, when traveling from $q$ to $p$. Here $m$ is the unique nonnegative integer such that

$$2^m d(q) \leq T < 2^{m+1} d(q).$$

Set $p = p_0$ and $q = q_0$. Thus we note that the curve $[p,q]_k$ has been divided into $(n+m+1)$ non-overlapping (modulo end points) subcurves:

$$[p_0,p_1]_k, \ldots, [p_{n-1},p_n]_k, \ [p_n,q_m]_k, \ [q_m,q_{m-1}]_k, \ldots, [q_1,q_0]_k.$$

Moreover, all subcurves in the above are also quasihyperbolic geodesics between their respective endpoints, and for any one of the above, denoted by $[x,y]_k$, we have

$$(3.5) \qquad d(u) \leq 2\min\{d(x), d(y)\} \text{ for all } u \in [x,y]_k.$$

Then for every $0 \leq i \leq n-1$, by (3.5), we may apply the assumption to the subarc $[p_i, p_{i+1}]_k$ and obtain by (2.2) that

$$(3.6) \qquad \log\left(1 + \frac{\ell([p_i, p_{i+1}]_k)}{d(p_i)}\right) \leq k(p_i, p_{i+1}) \leq A,$$

which ensures that

$$(3.7) \qquad \ell([p_i, p_{i+1}]_k) \leq e^A d(p_i).$$

Similarly, for each $0 \leq j \leq m-1$, we get

$$(3.8) \qquad \ell([q_{j+1}, q_j]_k) \leq e^A d(q_j) \quad \text{and} \quad \ell([p_n, q_m]_k) \leq e^A d(p_n).$$

We are now in a position to complete the proof of Theorem 1.2. For each $v \in [p,q]_k$, we have three cases: $v \in [p_i, p_{i+1}]_k$ or $[p_n, q_m]_k$ or $[q_{j+1}, q_j]_k$ for some $0 \leq i \leq n-1$ and $0 \leq j \leq m-1$. It is enough to consider the former two cases because the last one follows from a similar argument as the first case.

On the one hand, if $v \in [p_i, p_{i+1}]_k$ for some $0 \leq i \leq n-1$, then by (2.1) and (3.6) we have

$$\left|\log \frac{d(v)}{d(p_i)}\right| \leq k(p_i, v) \leq k(p_i, p_{i+1}) \leq A$$

and thus,

$$(3.9) \qquad d(p_i) \leq e^A d(v).$$

This, together with (3.4) and (3.7), shows that

$$\ell([p,v]_k) \leq \sum_{s=0}^{i} \ell([p_s, p_{s+1}]_k) \leq e^A \sum_{s=0}^{i} d(p_s) \leq 2e^A d(p_i) \leq 2e^{2A} d(v),$$

as desired.

On the other hand, if $v \in [p_n, q_m]_k$, then we easily see from a similar argument as (3.9) that $d(p_n) \leq e^A d(v)$. Using this and (3.8), we obtain

$$\begin{aligned} \ell([p,v]_k) &\leq \sum_{s=0}^{n} \ell([p_s, p_{s+1}]_k) + \ell([p_n, q_m]_k) \\ &\leq 2e^A d(p_n) + e^A d(p_n) \\ &\leq 3e^{2A} d(v). \end{aligned}$$

Therefore, we show that the quasihyperbolic geodesic $[p,q]_k$ is a double $M$-cone arc with $M = 3e^{2A}$. Hence Theorem 1.2 follows. $\square$

3.2. **Gromov hyperbolic John spaces are quasihyperbolic John.** In this part, we prove Theorem 1.3 with the aid of Theorem 1.2. Together with Lemma 3.1 and the stability property of quasigeodeiscs in Gromov hyperbolic spaces, we verify the criterion for quasihyperbolic John spaces.

**Lemma 3.2.** *Let $(X, d)$ be a minimally nice $a$-John Gromov $\delta$-hyperbolic space, and let $[x, y]_k$ be a quasihyperbolic geodesic in $X$ connecting $x$ and $y$. If $d(u) \leq 2 \min\{d(x), d(y)\}$ for all $u \in [x, y]_k$, then there is a positive number $A$ such that $k(x, y) \leq A$ with $A$ depending only on $a$ and $\delta$.*

*Proof.* First, take a double $a$-cone arc $\gamma$ connecting $x$ and $y$ in $X$ and let $x_0$ be the midpoint of $\gamma$. Then pick another two quasihyperbolic geodesics $[x, x_0]_k$ and $[y, x_0]_k$ joining $x_0$ to $x$ and $y$, respectively. Consider the quasihyperbolic geodesic triangle
$$\Delta = [x_0, x]_k \cup [x, y]_k \cup [y, x_0]_k.$$
There is a tripod map $f : \Delta \to T$, where $T$ is a tripod consisting of three line segments in $\mathbb{R}^2$ satisfying the following two properties: $f$ is an isometry on each of the sides of $\Delta$ and $f(u) = f(v)$ implies $k(u, v) \leq 4\delta$, see [2, Proposition 3.1]. Let $u \in [x, y]_k$, $v \in [x, x_0]_k$ and $w \in [x_0, y]_k$ be the points of the quasihyperbolic geodesic triangle $\Delta$ whose image under $f$ is the origin of $T$. Thus we have
$$(3.10) \qquad \max\{k(u, v), k(u, w)\} \leq 4\delta.$$
Without loss of generality, we may assume that $k(x, u) \geq k(u, y)$. Thus
$$(3.11) \qquad k(x, y) \leq 2k(x, u).$$
This is because for the other case that $k(x, u) \leq k(u, y)$, one may obtain the desired result by a symmetric argument.

Next, by Lemma 3.1, one immediately see that $\gamma[x, x_0]$ is a quasihyperbolic $(\lambda, \mu)$-quasigeodesic with $\lambda = 3a$ and $\mu = 3a \log(3a)$. Moreover, according to Theorem 2.1, there is a point $v_0 \in \gamma[x, x_0]$ and a number $R = R(\lambda, \mu, \delta) = R(a, \delta)$ such that
$$k(v_0, v) \leq R.$$
This, together with (2.1) and (3.10), guarantees that
$$(3.12) \qquad \left|\log \frac{d(v_0)}{d(u)}\right| \leq k(v_0, u) \leq k(v_0, v) + k(v, u) \leq R + 4\delta,$$
and so
$$d(v_0) \leq e^{R+4\delta} d(u) \leq 2e^{R+4\delta} d(x),$$
where the last inequality follows from the assumption.

Moreover, by (3.11) and (3.12), a second application of Lemma 3.1 gives
$$\begin{aligned} k(x, y) &\leq 2k(x, u) \leq 2k(x, v_0) + 2k(v_0, u) \\ &\leq 6a \log\left(1 + \frac{ad(v_0)}{d(x)}\right) + 2R + 8\delta \\ &\leq 6a \log\left(1 + 2ae^{R+4\delta}\right) + 2R + 8\delta := A \end{aligned}$$
and we complete the proof of Lemma 3.2. $\square$

**Proof of Theorem 1.3.** By Lemma 3.2 and Theorem 1.2, Theorem 1.3 follows immediately and one finds that every Gromov hyperbolic John space is quasihyperbolic John. □

3.3. **Quasihyperbolic John and GHS spaces.** In this subsection, we investigate the relationships between GHS spaces, linearly locally connected spaces and quasihyperbolic John spaces. Our goal is to prove Theorem 1.4. We start with several definitions.

**Definition 3.3.** Let $(X, d, \mu)$ be a metric measure space. We say that $X$ is *Ahlfors regular* if there are constants $Q > 1$ and $C > 0$ such that for each $x \in X$ and $0 < r \leq \mathrm{diam}(X)$, we have $C^{-1} r^Q \leq \mu(B(x,r)) \leq C r^Q$.

**Definition 3.4.** Let $(X, d)$ be a minimally nice space, and $C_{gh} \geq 1$ be a constant. We say that $X$ satisfies the $C_{gh}$-*Gehring-Hayman condition*, if for all $x$, $y$ in $X$ and for each quasihyperbolic geodesic $\gamma$ joining $x$ and $y$, we have $\ell(\gamma) \leq C_{gh} \ell(\beta_{x,y})$, where $\beta_{x,y}$ is a curve joining $x$ and $y$ in $X$.

**Definition 3.5.** Let $(X, d)$ be a minimally nice space, and $C_{bs} \geq 1$ be a constant. We say that $X$ satisfies the $C_{bs}$-*ball separation condition*, if for all $x$, $y$ in $X$, for each quasihyperbolic geodesic $\gamma$ joining $x$ and $y$, and for every $z \in \gamma$,
$$B(z, C_{bs} d(z)) \cap \beta_{x,y} \neq \emptyset,$$
where $\beta_{x,y}$ is any curve joining $x$ and $y$ in $X$.

Next, we introduce some facts about metric geometry of minimally nice spaces, see [3, 20]. Let $(X, d)$ be a minimally nice space. The length distance $\ell(x, y)$ of $(X, d)$ is defined by
$$\ell(x, y) = \inf\{\ell(\gamma)\},$$
where the infimum is taken over all rectifiable curves $\gamma$ in $X$ with endpoints $x$ and $y$. The space $(X, d)$ is called a *length* space if $d = \ell$.

Recall that a GHS space is a length space satisfying both the Gehring-Hayman and ball separation conditions defined as above. We first prove that a GHS space satisfying the $LLC_2$ condition is quasihyperbolic John.

**Proposition 3.6.** *Let $(X, d)$ be a minimally nice length space. If $(X, d)$ is GHS and $C_0$-$LLC_2$, then $X$ is quasihyperbolic John.*

*Proof.* For each pair of points $x, y \in X$ and for any quasihyperbolic geodesic $[x, y]_k = \alpha$ joining $x$ to $y$, we need to find an upper bound for the constant $C$, where $z \in \alpha$ and
$$\min\{\ell(\alpha[x,z]), \ell(\alpha[z,y])\} = C d(z).$$
Since $(X, d)$ satisfies the $C_{gh}$-Gehring-Hayman condition, we have
$$\ell(\alpha[x,z]) \leq C_{gh} d(x,z) \quad \text{and} \quad \ell(\alpha[y,z]) \leq C_{gh} d(y,z),$$
because both the subarcs $\alpha[x,z]$ and $\alpha[z,y]$ are quasihyperbolic geodesics. This in turn implies that
$$\min\{d(x,z), d(z,y)\} \geq \frac{C}{C_{gh}} d(z)$$

and so
$$x, y \in X \setminus \overline{B}\Big(z, \frac{C}{C_{gh}}d(z)\Big).$$

Now the $C_0$-$LLC_2$ property of $(X,d)$ ensures that there is a curve $\gamma$ joining $x$ to $y$ such that
$$\gamma \subset X \setminus \overline{B}\Big(z, \frac{C}{C_0 C_{gh}}d(z)\Big).$$

Moreover, since $\alpha$ is a quasihyperbolic geodesic and since $(X,d)$ satisfies the $C_{bs}$-ball separation condition, it follows that
$$\gamma \cap B(z, C_{bs}d(z)) \neq \emptyset.$$

This, together with the above fact, shows that
$$C \leq C_0 C_{gh} C_{bs},$$

as required. The proof of Proposition 3.6 is complete. □

We now derive a couple of corollaries of the above result. The first is an analog of [5, Theorem 4.2]. Recall that $(X,d)$ is *c-locally externally connected*, abbreviated *c*-LEC, provided $c \geq 1$ and the $LLC_2$ property holds for all points $x \in X$ and for all $r \in (0, d(x)/c)$, see [5].

**Corollary 3.7.** *A minimally nice space $(X,d)$ is uniform and LEC if and only if $(X,d)$ is quasiconvex and $(X,l)$ is $LLC_2$ and GHS, where $l$ denotes the length metric of $X$.*

*Proof.* Note that the sufficiency follows from Proposition 3.6. By [1, Theorem 2.3], we know that a length space satisfying a slice condition is GHS. Thus the necessity part follows from [5, Theorem 4.2]. □

The second is an analog of [3, Proposition 7.12] in the setting of metric spaces.

**Corollary 3.8.** *Let $(X, d, \mu)$ be a minimally nice, Ahlfors regular metric measure space. If $(X, l)$ is $LLC_2$ and if $(X, k)$ is a roughly starlike, $\delta$-hyperbolic space, then $(X, l)$ is uniform.*

*Proof.* We first note that $(X, l, \mu)$ is a minimally nice, upper $Q$-regular length space with $Q > 1$ and the measure $\mu$ is Ahlfors $Q$-regular on Whitney type balls in the sense of [24]. Thus applying [24, Theorem 5.1] to the space $(X, l, \mu)$, one finds that $(X, l)$ is a GHS space. Moreover, by Proposition 3.6, we know that every quasihyperbolic geodesic in $X$ is a double cone arc. This last fact and the fact that $(X, l)$ satisfies the Gehring-Hayman condition yield that $(X, l)$ is uniform, as required. □

Conversely, we have the following result which may be of independent interest.

**Proposition 3.9.** *Let $(X, d)$ be a minimally nice length space. If $X$ is quasihyperbolic a-John and b-LEC, then it is c-$LLC_2$ and satisfies the ball separation condition.*

*Proof.* We first check that $(X, d)$ is $c$-$LLC_2$ with $c = 3ab(b+1)$. Let $x_0 \in X$ and $y, z \in X \setminus \overline{B}(x_0, r)$. Pick a double $a$-cone arc $\alpha$ joining $y$ and $z$ in $X$. Then for all $x \in \alpha$, we have

(3.13) $$\min\{d(x, y), d(x, z)\} \leq \min\{\ell(\alpha[x, y]), \ell(\alpha[x, z])\} \leq ad(x).$$

We consider two cases. If $\alpha \cap S(x_0, br/c) = \emptyset$, then there is nothing to prove because $\alpha \subset X$. If $\alpha \cap S(x_0, br/c) \neq \emptyset$, then there are two points $y_1, z_1 \in \alpha \cap S(x_0, br/c)$ (perhaps $y_1 = z_1$) such that both $\alpha[y, y_1]$ and $\alpha[z_1, z]$ lie outside of $B(x_0, br/c)$. Thus, by (3.13) we get

$$ad(y_1) \geq \min\{d(y_1, y), d(y_1, z)\} \geq \min\{d(x_0, y), d(x_0, z)\} - d(x_0, y_1) \geq r - br/c,$$

which implies that

$$d(x_0) \geq d(y_1) - d(x_0, y_1) \geq (r - br/c)/a - br/c > 2b^2 r/c,$$

where the last inequality follows from our choice of the constant $c$. That is to say,

$$\overline{B}(x_0, br/c) \subset B\left(x_0, \frac{d(x_0)}{2b}\right).$$

Applying the $LEC$ condition, there is a curve $\beta \subset X \setminus \overline{B}(x_0, r/c)$ joining $y_1$ and $z_1$. Therefore, we get a curve

$$\gamma := \alpha[y, y_1] \cup \beta \cup \alpha[z_1, z],$$

which connects $y$ and $z$ in $X \setminus \overline{B}(x_0, r/c)$. This shows the first assertion.

We are thus left to prove that $(X, d)$ satisfies the ball separation condition. This can be seen as follows. For each pair of points $x, y \in X$ and for every quasihyperbolic geodesic $\alpha = [x, y]_k$ joining $x$ and $y$, we have

$$\{x, y\} \cap B(z, ad(z)) \neq \emptyset \quad \text{for all} \ \ z \in \alpha,$$

because $\alpha$ is a double $a$-cone arc. This deduces that

$$\min\{d(x, z), d(z, y)\} \leq \min\{\ell(\alpha[x, z]), \ell(\alpha[z, y])\} \leq ad(z),$$

which completes the proof. $\square$

Finally, we shall show that GHS John spaces are also quasihyperbolic John. By [1, Theorem 6.1], we see that a GHS space is Gromov hyperbolic with respect to its quasihyperbolic metric. From this fact and Theorem 1.3, it follows that GHS and John spaces are quasihyperbolic John. In the following, we give a direct proof for this conclusion.

**Proposition 3.10.** *Let $(X, d)$ be a minimally nice length space. If $(X, d)$ is GHS and $a$-John, then $X$ is quasihyperbolic John.*

*Proof.* Fix $x, y \in X$ and take a quasihyperbolic geodesic $\alpha = [x, y]_k$ connecting $x$ to $y$. Since $X$ is $a$-John, there is a double $a$-cone arc $\gamma$ joining $x$ and $y$. Because $(X, d)$ satisfies the $C_{bs}$-ball separation condition, for each $z \in \alpha$, there is $z_0 \in \gamma$ such that

(3.14) $$d(z, z_0) \leq C_{bs} d(z),$$

which implies
$$d(z_0) \leq d(z, z_0) + d(z) \leq (1 + C_{bs})d(z).$$

Moreover, since $\gamma$ is a double $a$-cone arc, thus we have
$$\min\{\ell(\gamma[x, z_0]), \ell(\gamma[z_0, y])\} \leq ad(z_0) \leq a(1 + C_{bs})d(z).$$

This combined with (3.14), ensures that

$$(3.15) \quad \begin{aligned} \min\{d(x,z), d(z,y)\} &\leq d(z, z_0) + \min\{\ell(\gamma[x, z_0]), \ell(\gamma[z_0, y])\} \\ &\leq (a+1)(1+C_{bs})d(z). \end{aligned}$$

On the other hand, since $(X, d)$ satisfies the $C_{gh}$-Gehring-Hayman inequality, it follows that
$$\ell(\alpha[x, z]) \leq C_{gh}d(x, z) \quad \text{and} \quad \ell(\alpha[y, z]) \leq C_{gh}d(y, z),$$

because both the subarcs $\alpha[x, z]$ and $\alpha[z, y]$ are quasihyperbolic geodesics. Therefore, by (3.15) we obtain
$$\min\{\ell(\alpha[x, z]), \ell(\alpha[z, y])\} \leq C_{gh}(a+1)(1+C_{bs})d(z).$$

Proposition 3.10 follows. $\square$

**Corollary 3.11.** *A minimally nice length space is uniform if and only if it is John and GHS.*

**Proof of Theorem 1.4.** Theorem 1.4 follows from Propositions 3.6 and 3.10. $\square$

## 4. Proofs of Corollaries 1.6 and 1.7

**4.1. Proof of Corollary 1.6.** Let $(X, d, \mu)$ be a minimally nice, Ahlfors regular metric measure space. To prove the necessity, we assume that $X$ is an $a$-uniform space. By definition, the quasiconvexity and John properties follow immediately. Thanks to [3, Theorem 3.6], $(X, k)$ is a geodesic Gromov hyperbolic space, where $k$ is the quasihyperbolic metric of $X$. It remains to show the roughly starlikeness of $(X, k)$. If $X$ is bounded, this again follows from [3, Theorem 3.6].

We next consider the case that $X$ is unbounded. Take a base point $p \in \partial X$. Let $(X, \widehat{d}_p)$ be the sphericalized space of $(X, d)$ with respect to $p$, for the definition see [7, Section 3.B]. According to [7, Theorem 5.5], we know that $(X, \widehat{d}_p)$ is a bounded $B$-uniform metric space with $B = B(a)$ because $(X, d)$ is unbounded. Arguing as the former case, $(X, \widehat{k}_p)$ is a geodesic Gromov hyperbolic metric space and $K'$-roughly starlike with respect to a point $w_0 \in X$, where $\widehat{k}_p$ is the quasihyperbolic metric of $X$ with respect to $\widehat{d}_p$.

For any $x \in X$, there is a $\xi$ on the Gromov boundary of hyperbolic space $(X, \widehat{k}_p)$ and a $\widehat{k}_p$-quasihyperbolic geodesic ray $[w_0, \xi]_{\widehat{k}_p}$ emanating from $w_0$ to $\xi$ such that
$$\widehat{k}_p(x, [w_0, \xi]_{\widehat{k}_p}) \leq K'.$$

Moreover, [7, Theorem 4.12] ensures that the identity map $(X, k) \to (X, \widehat{k}_p)$ is $(80a)$-bilipschitz. So we may assume that $\xi$ also belongs to the Gromov hyperbolic of the hyperbolic space $(X, k)$ and $[w_0, \xi]_{\widehat{k_p}}$ is a $k$-quasihyperbolic quasigeodesic ray. Then take another quasihyperbolic geodesic ray $[w_0, \xi]_k$ connecting $w_0$ and $\xi$. By the extended stability theorem (cf. [28, Theorem 6.32]) of Gromov hyperbolic spaces, there is a positive number $R = R(a)$ such that the Hausdorff distance between the two curves $[w_0, \xi]_k$ and $[w_0, \xi]_{\widehat{k_p}}$ with respect to the quasihyperbolic metric $k$ is bounded above by $R$. Therefore, we obtain

$$k(x, [w_0, \xi]_k) \leq 80AK' + R.$$

Hence we are done and the necessity follows.

For the sufficiency, we only need to show that every quasihyperbolic geodesic in $X$ is a uniform arc. On the one hand, Theorem 1.3 ensures that every quasihyperbolic geodesic is a double $M$-cone arc. On the other hand, according to [24, Theorem 5.1], the Gehring-Hayman condition follows. This, together with the quasiconvexity of $(X, d)$, guarantees that each quasihyperbolic geodesic is quasiconvex, as desired. □

4.2. **Proof of Corollary 1.7.** First, thanks to [1, Theorem 2.3], we know that $(X, l)$ satisfies the Gehring-Hayman and ball separation conditions, where $l$ is the length metric on $X$ induced by $d$. Thus $(X, l)$ is a GHS space. Now it follows from Proposition 3.10 that $X$ is quasihyperbolic John. Since $X$ is quasiconvex and since $(X, l)$ satisfies the Gehring-Hayman condition, we find that $X$ is uniform. Hence, Corollary 1.7 follows. □

QINGSHAN ZHOU, SCHOOL OF MATHEMATICS AND BIG DATA, FOSHAN UNIVERSITY, FOSHAN, GUANGDONG 528000, PEOPLE'S REPUBLIC OF CHINA
 *E-mail address*: `qszhou1989@163.com; q476308142@qq.com`

SAMINATHAN PONNUSAMY, DEPARTMENT OF MATHEMATICS, INDIAN INSTITUTE OF TECHNOLOGY MADRAS, CHENNAI 600036, INDIA
 *E-mail address*: `samy@iitm.ac.in`